\documentclass[12pt]{amsart}
\usepackage{amsmath,amsthm,amsfonts,amssymb,times}

\DeclareMathAlphabet{\curly}{U}{rsfs}{m}{n}

\theoremstyle{remark}

\theoremstyle{plain}

\newtheorem{lem}{Lemma}[section]
\newtheorem{thm}{Theorem}

\newtheorem{conj}{Conjecture}
\numberwithin{equation}{section}

\newcommand{\QQ}{{\mathbb Q}}
\newcommand{\RR}{{\mathbb R}}
\newcommand{\CC}{{\mathbb C}}

\newcommand{\rad}{{\rm rad}}

\newcommand{\be}{\begin{equation}}
\newcommand{\ee}{\end{equation}}
\newcommand{\benn}{\begin{equation*}}
\newcommand{\eenn}{\end{equation*}}
\newcommand{\bal}{\begin{align*}}
\newcommand{\ea}{\end{align*}}
\newcommand{\eal}{\ensuremath{\end{align*}}}
\newcommand{\bea}{\begin{eqnarray}}
\newcommand{\eea}{\end{eqnarray}}

\newcommand{\lam}{\ensuremath{\lambda}}

\renewcommand{\a}{\ensuremath{\alpha}}

\newcommand{\del}{\ensuremath{\delta}}
\newcommand{\eps}{\ensuremath{\varepsilon}}
\renewcommand{\(}{\left(}
\renewcommand{\)}{\right)}

\newcommand{\pfrac}[2]{\left(\frac{#1}{#2}\right)}

\def\RR{{\mathcal R}}
\def\SS{{\mathcal S}}
\def\UU{{\mathcal U}}
\def\TT{{\mathcal T}}
\def\PP{{\mathcal P}}
\def\EE{{\mathcal E}}
\def\CC{{\mathcal C}}
\def\MM{{\mathcal M}}
\def\QQ{{\mathcal Q}}

\renewcommand{\le}{\leqslant}

\renewcommand{\ge}{\geqslant}

\begin{document}

\title{Common values of the arithmetic  functions $\phi$ and $\sigma$}
\author{Kevin Ford}
\address{\noindent Department of Mathematics, 1409 W. Green St., University
of Illinois at Urbana-Champaign, Urbana, IL 61801, USA}
\email{ford@math.uiuc.edu}
\author{Florian Luca}
\address{Instituto de Matem{\'a}ticas,
Universidad Nacional Autonoma de M{\'e}xico,
C.P. 58089, Morelia, Michoac{\'a}n, M{\'e}xico} 
\email{fluca@matmor.unam.mx}
\author{Carl Pomerance}
\address{Department of Mathematics, Dartmouth college, Hanover, NH 03755, USA}
\email{carl.pomerance@dartmouth.edu}
\date{\today}
\begin{abstract}  We show that the equation $\phi(a)=\sigma(b)$ has infinitely 
many solutions, where $\phi$ is Euler's totient function and $\sigma$ is the
sum-of-divisors function.  This proves a 50-year old conjecture of Erd\H os.
Moreover, we show that for some $c>0$, there are infinitely many
integers $n$ such that
$\phi(a)=n$ and $\sigma(b)=n$ each have more than $n^c$
solutions.  The proofs rely on the recent work of
the first two authors and Konyagin on the distribution of
primes $p$ for which a given prime divides some iterate of $\phi$ at $p$,
and on a result of Heath-Brown
connecting the possible existence of Siegel zeros with the distribution
of twin primes.
\end{abstract}

\subjclass[2000]{11A25, 11N25, 11N64}

\thanks{The research of K.~Ford was supported in part
by NSF grant DMS-0555367, 
that of F.~Luca was supported in part
by projects PAPIIT 100508 and SEP-CONACyT 79685, and that of
C.~Pomerance was supported in part by NSF grant DMS-0703850.
This paper was begun while the first two authors were visiting
Dartmouth College; they thank the Mathematics Department for
its hospitality.  The authors thank Paul Pollack for some helpful
conversations.}

\maketitle


\section{Introduction}\label{sec:intro}

Two of the oldest and most studied functions in the theory of numbers are the
sum-of-divisors function $\sigma$ and Euler's totient function $\phi$.
Over 50 years ago, Paul Erd\H os conjectured that the ranges  of $\phi$ and 
$\sigma$ have an infinite intersection
(\cite[p.~172]{E59}, \cite[p.~198]{SS}).   This conjecture follows easily from
some famous unsolved problems.  For example,
if there are infinitely many pairs of twin primes $p$, $p+2$,
 then $\phi(p+2)=p+1=\sigma(p)$, and if there are infinitely many
Mersenne primes $2^p-1$, then $\sigma(2^p-1)=2^p=\phi(2^{p+1})$.
Results from \cite{F98} indicate that typical values taken by $\phi$ and by $\sigma$ have a similar multiplicative structure; hence, common values should be plentiful.
A short calculation reveals that there are 95145 common values of $\phi$
and $\sigma$ between 1 and $10^6$.  This is to be compared with a total
of 180184 $\phi$-values and 189511 $\sigma$-values in the same interval.
In \cite{EG}, the authors write that ``it is very annoying that we
cannot show that $\phi(a)=\sigma(b)$ has 
infinitely many solutions$\dots$ ."
Annoying of course, since it is so obviously correct!
Erd\H os knew (see \cite[sec. B38]{Guy}) that $\phi(a)=k!$ is solvable for every positive
integer $k$, so all one would have to do is show that $\sigma(b)=k!$
is solvable for infinitely many choices for $k$.  In fact, this
equation seems to be solvable for every $k\ne2$, but proving it seems
difficult.

The heart of the problem is to understand well 
the multiplicative structure of shifted primes $p-1$ and $p+1$.

In this note, we give an unconditional proof of the Erd\H os conjecture.  
Key ingredients in the proof are a very recent bound
on counts of prime chains from \cite{FKL} (see \S \ref{sec:chains}
for a definition) and estimates for
primes in arithmetic progressions.  The possible existence of
Siegel zeros (see \S \ref{sec:progressions} for a definition)
creates a major obstacle for the success of our argument.
Fortunately, Heath-Brown \cite{HB} showed that if
Siegel zeros exist, then there are infinitely many 
pairs of twin primes.
However, despite the influence of possible Siegel zeros,
our methods are completely effective.

\begin{thm}\label{thm1}
The equation $\phi(a)=\sigma(b)$ has infinitely many solutions.
Moreover, for some positive $\alpha$ and all large $x$, there are at least 
$\exp \left((\log\log x)^{\alpha}\right)$ integers
$n\le x$ which are common values of $\phi$ and $\sigma$.
\end{thm}

We also show that there are infinitely many integers $n$ which are
common values of $\phi$ and $\sigma$ in many ways.
Let $A(n)$ be the number of solutions of $\phi(x)=n$,
and let $B(n)$ be the number of solutions of $\sigma(x)=n$.  Pillai \cite{Pi}
showed in 1929 that the function $A(n)$ is unbounded, and in 1935, Erd\H os 
\cite{E35} showed that the inequality $A(n)>n^c$ holds infinitely often for some positive constant $c$.  
The proofs give analogous results for $B(n)$.  Numerical values of $c$
have been given by a number of people
(\cite{Ba}, \cite{Fr}, \cite{Po} and \cite{Wo}), the largest so far being
$c=0.7039$ which is due to Baker and Harman \cite{BH}.  The key to these
results is to show that there are many primes $p$ 
for which $p-1$ has only small prime factors.
Erd\H os \cite{E45} conjectured that for any constant $c<1$ the inequality 
$A(n)>n^c$ holds infinitely often. 

\begin{thm}\label{thm2}
For some 
positive constant $c$ there are infinitely many $n$ such that both inequalities $A(n)>n^c$ and $B(n)>n^c$ hold.  Moreover,
for some constant $a>0$, there are at least
$(\log\log x)^a$ such numbers $n\le x$, for all large $x$.
\end{thm} 

Necessary results on the distribution of primes in progressions, twin primes,
and prime chains are given in Sections \ref{sec:progressions}
and \ref{sec:chains}.
In Section \ref{sec:chains}, we prove Theorem~\ref{thm1}.
In Section \ref{sec:popular}, we present
 the additional arguments needed 
to deduce the conclusion of Theorem \ref{thm2}.  
Theorem \ref{thm2} resolves another conjecture of Erd\H os
(stated as Conjecture $C_8$ in \cite[p.~193]{SS}):
for each number $k$, there is some number $n$ with 
$A(n)>k$ and $B(n)>k$.  
Later, in Section \ref{sec:further}, we pose some additional
problems concerning common values of $\phi$ and $\sigma$.

We consider $n=\sigma\left(\prod_{p\in \SS} p\right) =\prod_{p\in \SS} (p+1)$,
where $\SS$ is a set of primes $p\le x$ for which all prime factors of
$p+1$ are small,
say $\le z$.  In this way, $n$ should be the product of some of the primes $\le z$, 
each to a possibly large power.  We deduce that $n$ is in the range of $\phi$
by exploiting the general implication
\be\label{rad}
\phi(\rad(m))\mid m \, \implies \, m = \phi\pfrac{m{\cdot}\rad(m)}{\phi(\rad(m))},
\ee
where $\rad(m)$ is the product of the distinct prime factors of $m$.
Let $v_q(m)$ denote the exponent of $q$ in the factorization of $m$.
We expect for $n=\sigma\left(\prod_{p\in\SS}p\right)$
that $v_q(\phi(\rad(n)))\le v_q(n)$ for $q\le z$;
hence, the hypothesis in \eqref{rad} should hold. 
Turning this into a proof requires lower bounds of the 
expected order for the number of $p\in \SS$
for which $q\mid p+1$.

We remark that by our proofs below, the numbers $n$ which are
constructed for Theorems \ref{thm1} and \ref{thm2} are also values 
taken by the Carmichael function $\lambda(m)$, the largest order of
an element of $(\mathbb{Z}/m\mathbb{Z})^*$.  Moreover, for the $n$ in
Theorem \ref{thm2}, there are at least $n^c$ such values $m$.
We thank Bill Banks for this observation. 

%
%
\section{Primes in progressions}\label{sec:progressions}
%
%

Throughout, constants implied by $O$, $\ll$, $\gg$, and $\asymp$
notation are absolute
unless otherwise noted.  Bounds for implied constants, as well as 
positive quantities introduced later, are effectively computable.
Symbols $p,q,r$ always denote primes,
and $P(m)$ is the largest prime factor of an integer $m>1$.
Let $\pi(x;m,a)$ be the number of primes $p\le x$ with $p\equiv
a\pmod{m}$, and let
$$
\psi(x;m,a) = \sum_{\substack{n\le x \\ n\equiv
 a\!\!\!\!\pmod{m}}} \Lambda(n),
$$
where $\Lambda$ is the von Mangoldt function.
The behavior of $\pi(x;m,a)$ and $\psi(x;m,a)$ are intimately connected to the
distribution of zeros of Dirichlet $L$-functions.  Of particular importance are
possible zeros near the point $1$.  
Let $\CC(m)$ denote the set of primitive characters modulo $m$.
It is known (cf.~\cite[Ch.~14]{Da}) that 
for some constant $c_0>0$ and every $m\ge 3$, there is at most one 
zero of $\prod_{\chi\in C(m)} L(s,\chi)$ in the region
\be
\label{zf}
\Re s \ge 1 - \frac{c_0}{\log (m(|\Im s|+1))}.
\ee
Furthermore, if this ``exceptional zero''
$\beta$ exists, it is real, it is a zero of $L(s,\chi)$ for a real character $\chi\in \CC(m)$,
and
\be\label{beta}
\beta \le 1-\frac{c_1}{m^{1/2}\log^2 m}
\ee
for some positive constant $c_1$.
Better upper bounds on $\beta$ are known (Siegel's Theorem, \cite[Ch. 21]{Da}),
but these are ineffective.
The ``exceptional moduli'' $m$, for which an exceptional $\beta$ exists,
must be quite 
sparse, as the following classical results show (\cite[Ch.~14]{Da}).

\begin{lem}[Landau]\label{Landau}
For some constant $c_2>0$, if $3\le m_1 < m_2$, $\chi_1\in \CC(m_1)$ and 
$\chi_2\in{\CC(m_2)}$, then there is at most one zero $\beta$ of 
$L(s,\chi_1) L(s,\chi_2)$ with $\beta > 1-{c_2}/{\log(m_1 m_2)}$.
\end{lem}

We immediately obtain

\begin{lem}[Page]\label{Page}
For any $M\ge 3$,
$$
\prod_{m\le M} \prod_{\chi\in \CC(m)} L(s,\chi)
$$
has at most one zero in the interval $[1-{(c_2/2)}/{\log M},1]$.
\end{lem}

It is known after
McCurley \cite{Mc} that $c_0=1/9.645908801$ holds in \eqref{zf},
while Kadiri \cite{Ka} has shown we may take
$c_0=1/6.397$, and in Lemmas~\ref{Landau}, \ref{Page} 
we may take $c_2=1/2.0452$.

The Riemann hypothesis for Dirichlet $L$-functions implies that no exceptional
zeros can exist.  If there is an infinite sequence of 
integers $m$ and associated zeros $\beta$ satisfying $(1-\beta)\log m \to 0$,
such zeros are known as Siegel zeros, and their existence would have 
profound implications on the distribution of primes in arithmetic progressions 
(\cite[(9) in Ch. 20]{Da}).
As mentioned before, Heath-Brown showed that the existence of Siegel zeros
implies that there are infinitely many prime twins.

\begin{lem}[{\cite[Corollary 2]{HB}}]\label{HB}
If $\chi\in \CC(m)$ and $L(\beta,\chi)=0$ for $\beta = 1 - \lam (\log m)^{-1}$, then
for $m^{300} < z \le m^{500}$, the number of primes $p\le z$ with $p+2$ prime
is
$$
C \frac{z}{\log^2 z} + O\( \frac{\lam z}{\log^2 z} \),
\quad\hbox{ where }C=2\prod_{p>2} (1-(p-1)^{-2})=1.32\ldots\,.
$$
\end{lem}

If Siegel zeros do not exist, there still may be some 
Dirichlet $L$-function zeros with real part $>1/2$, which
would create irregularities in the distribution of primes in some
progressions.  Such progressions, however, would have moduli
larger than a small power of $x$.  We state here a character sum 
version of this result, due to Gallagher (see the proof of \cite[Theorem 7]{G}).  Let
$$
\psi(x,\chi)= \sum_{n\le x} \Lambda(n) \chi(n), \qquad  {\text{\rm and}}\qquad
\Psi(x,m) = \sum_{\chi\in \CC(m)} |\psi(x,\chi)|.
$$

\begin{lem}\label{Gallagher}
If $c_2$ is as in Lemma \ref{Landau}, then for
every $\lam \in (0, c_2/2]$ and $\eps>0$, there are constants $1\ge\a>0$ 
and $x_0$ so that for $x\ge x_0$,
$$
\sum_{\substack{3\le m\le x^{\a} \\ m\ne m_0}} \Psi(x,m) \le \eps x.
$$
Here $m_0$ corresponds to the conductor of a Dirichlet character $\chi$ 
for which $L(\beta,\chi)=0$ for some $\beta>1-\lam/\log(x^{\a})$.  
If there is no such zero, set $m_0=0$.
\end{lem}

We remark that $m_0$, if it exists, is unique by Lemma \ref{Page}.

We also know that $\Psi(x,m)$ is small for most $m\in (x^\a,x^{1/2-\del}]$
if $\del>0$ is fixed.  This follows from the next lemma which is a key ingredient in
the proof of the Bombieri-Vinogradov theorem.

\begin{lem}\label{Psi}
For $1\le M\le x$,
$$
\sum_{m\le M} \Psi(x,m) \ll \( x + x^{5/6}M +
x^{1/2}M^2\)\log^4 x.
$$
\end{lem}

\begin{proof}
This is \cite[Ch.~28, (2)]{Da}.
\end{proof}

For positive reals $\delta,\gamma,y,x$, with $1\le y\le x^{1/2-\delta}$,
and a nonzero integer $a$, define
\begin{align*}
S_q(x;\del,a)&=\# \{ p\le x : P(p+a) \le x^{1/2-\del}, q\mid p+a \}, \\
{\EE}(x,y;\del,\gamma) &= \left\{ q\le y : 
S_q(x;\del,1) \le \frac{\gamma x}{q\log x} 
\text{ or }
S_q(x;\del,-1) \le \frac{\gamma x}{q\log x} 
\right\}.
\end{align*}

We say that a real number $x$ is $(\a,\eps)$-good if 
$\Psi(x;m) \le \eps x$ for $3\le m\le x^{\a}$.  
Roughly speaking, this means that the exceptional modulus
in Lemma \ref{Gallagher} doesn't exist (for appropriate $\lam$).

\begin{lem}\label{E}
There are absolute constants $\del>0$ and $\gamma>0$ so that the following
holds.  For every $\a>0$, there are constants $\eta>0$
and $x_1>0$ so that if $x\ge x_1$ and $x$ is $(\a,\frac1{10})$-good, then
for all $y\le x^{1/2-\del}$,
$$
\# {\mathcal E}(x,y;\del,\gamma) \le y x^{-\eta} .
$$
\end{lem}

\begin{proof} We may assume that $0<\del < 1/6$.
Let $k$ be a positive integer such that
$Q=2^{-k} x^{1/2-\del}\ge 1$. Let $R_1=\max\{Q^{-1} x^{1/2-5\del/4},\,x^{\del/4}\}$,
and let $R_2=R_1 x^{\del/4}$.
By standard estimates (\cite[(3) in Ch. 20]{Da}),
if $q\in(Q,2Q]$ and $r\in(R_1,R_2]$, then for $a=\pm1$,
\be\label{psi}
\left| \psi(x;qr,a) - \frac{x}{\phi(qr)} \right| \le
\frac1{\phi(qr)}\left({\Psi(x,q)+\Psi(x,r)+\Psi(x,qr)+O(x/\log x)}\right).
\ee
Let ${\mathcal E}_1(Q)=\{q\in (Q,2Q]:\Psi(x,q) > x/10\}$.  
Since $x$ is $(\a,\frac1{10})$-good, we have ${\mathcal E}_1(Q)=\emptyset$ when
$Q\le \frac12 x^{\a}$.  Otherwise, by Lemma \ref{Psi},
$$
\#{\mathcal E}_1(Q) \ll \(1+Qx^{-1/6} + Q^2 x^{-1/2} \)\log^4 x \ll Q(x^{-\del}+x^{-\a})\log^4 x.
$$
Let ${\mathcal E}_2(Q)=\{q\in(Q,2Q]: \Psi(x,qr)> x/10
\text{ for at least }R_1 x^{-\del/8}\text{ primes }r\in (R_1,R_2]\}$.
By Lemma \ref{Psi} and the inequality $R_2 Q \le x^{1/2-\del/2}$,
$$
\#{\mathcal E}_2(Q) \ll \frac{(x+x^{5/6}R_2 Q + x^{1/2} (R_2 Q)^2)\log^4 x}
{R_1 x^{1-\del/8}} \ll Q x^{-\del/8} \log^4 x.
$$
Also, by Lemma \ref{Psi},
$$
\# \{r\in( R_1,R_2] : \Psi(x,r)\ge x/10 \} \ll \( 1 + x^{-1/6}R_2 + x^{-1/2} R_2^2 \)
\log^4 x \ll R_1x^{-\del/2} \log^4 x.
$$
For each $q\in(Q,2Q]$ with 
$q\not \in {\mathcal E}_1(Q) \cup {\mathcal E}_2(Q)$, 
let 
\[
\RR(q)=\{r\in (R_1,R_2]: \Psi(x,qr)\le x/10,\, \Psi(x,r) \le x/10\}.
\]
By \eqref{psi}, for $r\in \RR(q)$ and $a=\pm1$,
\be
\label{piest}
\pi(x;qr,a) \ge \frac{\psi(x;qr,a)-O(\sqrt{x})}{\log x} \ge \frac{x}{2qr\log x}.
\ee
Also, by the above estimates and Mertens' formula,
\be\label{sum1r}
\sum_{r\in \RR(q)} \frac{1}{r} \ge \sum_{R_1 < r \le R_2}
 \frac{1}{r} - O(x^{-\del/8}\log^4 x) \ge \frac{\del}{2}.
\ee
Since $R_1\ge x^{\del/4}$, a shifted prime $p+a$ is 
divisible by at most $\lfloor\frac{4}{\del}\rfloor$ 
primes in $\RR(q)$.  Hence,
\begin{align*}
S_q(x;\del,a) &\ge\frac{\del}{4}\sum_{r\in\RR(q)} 
\left( \pi(x;qr,-a)-\#\UU(q,r) \right),\qquad \text{ where}\\
\UU(q,r) &= \{p\le x: qr | p+a, P(p+a)>x^{1/2-\del} \}.
\end{align*}
Since $r\le R_2 \le x^{1/2-\del}$,
if $p\in \UU(q,r)$, then $p+a=qrsb$, where $s>x^{1/2-\del}$ is
prime and 
$$
b\le \frac{x+1}{qrs} \le \frac{x+1}{x^{1-9\del/4}} \le x^{3\del}.
$$
For fixed $b,q,r,a$, we estimate the number of possible choices for 
$s$ using the sieve (\cite{HR}, Theorem 3.12).  We get
$$
\#\UU(q,r) \ll \sum_{b\le x^{3\del}} \frac{x}{bqr \log^2 ({x}/{bqr})}
\frac{b}{\phi(b)}
\ll \frac{x}{qr\log^2 x} \sum_{b\le x^{3\del}} \frac{1}{\phi(b)}
\ll \frac{\del x}{qr\log x}.
$$
For small enough $\del$, we then have $\#\UU(q,r) \le \frac{x}{4qr\log x}$, 
and we conclude from \eqref{piest}, \eqref{sum1r} that
$$
S_q(x;\del,a) \ge \frac{\del x}{16q\log x} \sum_{r\in \RR(q)} \frac{1}{r}
\ge \frac{\del^2 x}{32q\log x}.
$$
Finally, $\#{\mathcal E_1}(Q)+\#{\mathcal E}_2(Q) \le \frac14 Q x^{-\eta}$
for $\eta=\min\{{\a}/{2},{\del}/{9}\}$ and large $x$.  Summing over
choices of the dyadic interval $(Q,2Q]$ with
$Q\le y$ and $a\in\{-1,1\}$ finishes the proof.
\end{proof}

%
%
\section{Prime chains and the proof of Theorem \ref{thm1}}\label{sec:chains}
%
%

Suppose that $n$ is a positive integer with $\phi(\rad(n))\mid n$
and that $q$ is a prime with $q\nmid n$.  Then $n$ is not divisible
by any prime $t\equiv1\pmod q$, since otherwise $q\mid\phi(\rad(n))$,
which would imply that $q\mid n$.  Iterating, $n$ is not divisible
by any prime $t'\equiv1\pmod t$, where $t$ is a prime with $t\equiv1\pmod q$.
And so on.
Thus, the single nondivisibility assumption that $q\nmid n$, plus
the assumption that $\phi(\rad(n))\mid n$, forces any prime $t$
in any {\em prime chain} for $q$ to also not divide $n$.  We define
a prime chain as a sequence of primes $q=t_0,t_1,t_2,\dots$, where
each $t_{j+1}\equiv 1\pmod{t_j}$.  Alternatively, if $\phi_j$ is the
$j$-fold iterate of $\phi$, then a prime $t$ is in a prime chain for $q$
if $t=q$ or $q\mid\phi_j(t)$ for some $j$.

Let $\TT(y,q)$ be the set of primes $t\le y$ which are in a prime chain for $q$.
Crucial to our proof is the following estimate.

\begin{lem}[{\cite[Theorem 5]{FKL}}]\label{FKL}
For every $\eps>0$ there is a constant $C(\eps)$ 
so that if $q$ is prime and $y>q$, then $\#\TT(y,q) \le C(\eps) (y/q)^{1+\eps}$.
\end{lem}

\noindent
More estimates for counts of prime chains with various properties may be found in
\cite{B, EGPS, FKL, LP}.

%
%
\bigskip

We now proceed to prove Theorem \ref{thm1}.
There is an absolute constant $\lam_0>0$ so that if $\lam\le \lam_0$, then
the error term in
the conclusion of Lemma \ref{HB} is at most
$0.1{z}/{\log^2 z}$ in absolute value.
Let $\a>0$ and $x_0$ be the constants from Lemma \ref{Gallagher}
corresponding to $\eps=\frac{1}{10}$ and $\lam=\lam_0$, and let
$\del,\gamma,\eta$ and $x_1$ be the constants from Lemma \ref{E}.

Suppose $x\ge \max(x_0,x_1)$.
We show that there are many common values of $\phi$ and $\sigma$ which are
$\le e^{2x}$ by considering two cases.
First, suppose that $x$ is not $(\a,\eps)$-good.
Then for some $m\le x^\a$ and $\chi\in \CC(m)$,
$L(\beta,\chi)=0$ for some
$\beta \ge 1 - {\lam_0}/{\log(x^{\a})}$.  By \eqref{beta}, 
$$
m \gg \frac{\log^2 x}{(\log\log x)^4}.
$$
Let $z=m^{500}$.
By Lemma \ref{HB}, the set $\TT$ of primes $p\le z-1$ for which $p+2$
is also prime satisfies $\# \TT \ge {1.2 z}/{\log^2 z}$.
Let $1/500>\theta>0$ be a sufficiently small
constant and $x$ large depending on $\theta$.  If $q=P(p+1)\le z^\theta$,
then $p+1=qb$ where $b$ is free of prime factors in $(q,z^{1/4}]$.
The number of such $p\in\TT$ is, by an application of the 
large sieve \cite[p.~159]{Da},
$$
\ll \sum_{q\le z^{\theta}} \frac{z}{\log^3 z} \frac{\log q}{q} \ll
 \frac{\theta z}{\log^2 z}.
$$
Let $\SS=\{p\in \TT : P(p+1) > z^{\theta} \}$.  Choose $\theta$ so small 
that $\# \SS \ge {z}/{\log^2 z}$.
For $p\in\SS$, we have
$$
\#\{p'\in\SS:P(p+1)\mid p'+1\}\le\frac{z}{P(p+1)}<z^{1-\theta}.
$$
Hence, there is a set $\PP$ of primes in $\SS$ with 
$\#\PP=\lfloor{z^\theta}/{\log^2 z}\rfloor$, and such that for each $p\in\PP$, $P(p+1)\nmid p'+1$ for all $p'\in\PP$ different from $p$.
For any subset $\MM$ of $\PP$, let $n(\MM)=\prod_{p\in \MM}(p+1)$,
so that $n(\MM)=\sigma(\prod_{p\in\MM}p)=\phi(\prod_{p\in\MM} (p+2))$.
Furthermore,
since each factor $p+1$ in the product $n(\MM)$ has the unique
``marker prime" $P(p+1)$ which divides no other $p'+1$ in the product,
the numbers $n(\MM)$ are distinct as $\MM$ varies.  Since 
$n(\MM)\le z^{\#\PP} \le x^{500\a x^{500\theta\a}}\le e^x$ for $x$ large,  
there are at least
$2^{\# \PP} > \exp\{z^{\theta/2}\}$ common values of $\phi$ and $\sigma$ 
which are $\le e^{x}$.  Observing that $z>(\log x)^{999}$ completes the 
proof in this case.

Now assume that $x$ is $(\a,\eps)$-good.
Let $\EE=\EE(x,x^{1/2-\del};\del,\gamma)$ and let
\be\label{Tdef}
\TT = \bigcup_{q\in \EE} \TT(x^{1/2-\del},q).
\ee
Put
\be\label{S}
\SS = \{ p\le x : P(p+1)\le x^{1/2-\del}\text{ and }t\nmid p+1 
\text{ for all } t\in \TT\}.
\ee
By partial summation and Lemmas \ref{E}, \ref{FKL}, we have for each $\eps>0$,
$$
\sum_{t\in\TT}\frac1t\le\sum_{q\in\EE}\sum_{t\in\TT(x^{1/2-\del},q)}\frac1t
\ll_\eps\sum_{q\in\EE}\frac{x^{(1/2-\del)\eps}}{q^{1+\eps}}\\
\ll_\eps x^{(1/2-\del)\eps-\eta}.
$$
Thus, if $\eps$ is small enough and $x$ large, we have 
\be
\label{sum1/t}
\sum_{t\in\TT}\frac1t<\frac{\gamma}{20\log x}.
\ee
Using Lemma \ref{E}, $2\not\in\EE$, so that
$\#\{p\le x:P(p+1)\le x^{1/2-\del}\}>(\gamma/2)x/\log x$.  
Thus, 
\be\label{Sb}
\#\SS>\frac{\gamma x}{2\log x}-\sum_{t\in\TT}\frac{x}{t}\ge\frac{\gamma x}{3\log x}.
\ee

Let $p_j$ be the  $j$-th largest prime in $\SS$, and
$$
n_j  =  \sigma \Bigl( \prod_{p\in \SS-\{p_j\}} p \Big) =  \prod_{p\in
 \SS-\{p_j\}} (p+1).
$$
Clearly $B(n_j)\ge 1$.  Note that the prime factors of $n_j$ 
are $\le x^{1/2-\del}$, so that
\[
\phi(\rad(n_j)) \mid u!,
\]
where $u=\lfloor x^{1/2-\del}\rfloor$.
If $q\le x^{1/2-\delta}$ and
$q\in \TT$, then $q\nmid \phi(\rad(n_j))$.
If $q\not\in \TT$, we have
\be
\label{vq}
v_q(\phi(\rad(n_j))) \le v_q(u!)  \le \frac{x^{1/2-\del}}{q-1}.
\ee
On the other hand, for such $q$, Lemma \ref{E} and \eqref{sum1/t} imply
\be\label{vqn}
v_q(n_j) \ge \# \{ p\in \SS-\{p_j\}: q\mid p+1 \}  \ge 
\frac{\gamma x}{q\log x} - 1
- \sum_{t\in T} \frac{x}{qt}  \ge \frac{\gamma x}{2q\log x} 
\ee
for $x$ sufficiently large.
Therefore, comparing \eqref{vq} with \eqref{vqn}
we see that \eqref{rad} holds with $m=n_j$ and so $A(n_j)\ge 1$.
By the prime number theorem,
$n_j\le \prod_{p\le x} (p+1) \le e^{2x}$ if $x$ is large. 
The numbers $n_j$ are distinct,  hence there
are at least $\# \SS \ge (\gamma/3) x/\log x$ 
common values of $\phi$ and $\sigma$ less than
 $e^{2x}$.  This completes the proof of Theorem \ref{thm1}.

%
%
%
\section{Popular common values of $\phi$ and
$\sigma$}\label{sec:popular}
%
%

In this section, we combine the proof of Theorem \ref{thm1} with
a method of Erd\H os \cite{E35}.  A key estimate 
is \cite[Lemma 2]{E35}:
\be\label{smooth}
\# \{ n\le x : P(n) \le \log x \}  = x^{o(1)} \qquad (x\to \infty).
\ee
More results about the distribution of integers $n$ with $P(n)$ small
may be found in \cite{HT}.

Define $\lam=\lam_0$, $\a$, $x_0$, $x_1$ and $\eta$ as in
 the proof of Theorem \ref{thm1}.  Without loss of generality, 
suppose $\a \le \frac1{500}$.
Theorem \ref{thm2} is proved by considering the two cases, $x$ is not 
$(\a,\frac{1}{10})$-good and $x$ is $(\a,\frac{1}{10})$-good.
The next lemmas provide the necessary arguments.

\begin{lem}\label{commonHB}
For some absolute constants $c>0$ and $a>0$, 
if $0 < \a \le \frac{1}{500}$, $x$ is large (depending on $\a$) 
and not $(\a,\frac1{10})$-good then there are
at least $(\log x)^{a}$ 
integers $n\le e^x$ for which both $A(n)>n^c$ and $B(n)>n^c$.
\end{lem}

\begin{proof}
As in the proof of Theorem \ref{thm1}, by \eqref{beta} 
there is an exceptional modulus $m$ satisfying 
$$
\frac{\log^2 x}{(\log\log x)^4} \ll m \le x^{\a}
$$
and so that
\be\label{twins}
\# \{ p\le z : p+2\text{ prime} \} \ge \frac{z}{\log^2 z}, 
\quad z=m^{500}.
\ee
Let $\del$ be a positive, absolute constant.
Let $\mathcal{P}$ be the 
set of primes $p\le z$ with $p+2$ prime and $P(p+1)\le z^{1-\del}$.
If $p$ and $p+2$ are both prime
 and $P(p+1)>z^{1-\del}$, then
$p+1=qb$ for some prime $q$ and some $b\le z^{\del}$.
 By sieve methods (\cite[Theorem 2.4]{HR}),
for small enough $\del$, we have
\begin{align*}
\# \PP &\ge \frac{z}{\log^2 z} - \sum_{b\le z^{\del}} \#\{q\le z/b:q, ~
qb-1,~qb+1 \text{ prime} \}\\
&\ge \frac{z}{\log^2 z}  - O\left(  \sum_{b\le z^{\del}} \frac{z}{b\log^3 z}
\pfrac{b}{\phi(b)}^2 \right) 
\ge \frac{z}{\log^2 z} - O\pfrac{\del z}{\log^2 z} \ge \frac{z}{2\log^2 z}.
\end{align*}
Let $H=\lfloor z^{1-\del/2} \rfloor$ and $J=\lfloor\#\PP/H \rfloor $.
Define sets $\PP_j$, $1\le j\le J$, as follows: $\PP_1$ is the set of the
smallest $H$ primes in $\PP$, $\PP_2$ is the set of the next $H$ smallest
primes from $\PP$, etc.  Let $K=\lceil z^{1-\del}/\log z \rceil$.
We may assume that $x$ is large enough that $K\ge2$, so that if ${\mathcal M}$
is a set of $K$ primes from some $\PP_j$, then
\be
\label{calM}
n({\mathcal M})=
\sigma\left(\prod_{p\in{\mathcal M}}p\right)
=\phi\left(\prod_{p\in{\mathcal M}}(p+2)\right)
\le z^K,
\quad P(n({\mathcal M}))\le z^{1-\del}\le\log\left(z^K\right).
\ee
By \eqref{smooth}, the function $n(\cdot)$ maps sets ${\mathcal M}$
into a set of integers of cardinality $\le z^{\del K/6}$.  But the number of
$K$-element subsets $\MM$ of some $\PP_j$ is
$$
\binom{H}{K}\ge\left(\frac{H}{K}\right)^K\ge z^{\del K/2}
$$
for $x$ large.
Thus, for each $j\le J$ there is some $n_j$ such that $n_j=n({\mathcal M})$
for at least $z^{\del K/3}$
$K$-element subsets ${\mathcal M}$ of $\PP_j$.
We conclude from \eqref{calM} that both 
$A(n_j),B(n_j)\ge z^{\del K/3}\ge n_j^{\del/3}$.
Since $n_1<n_2<\dots<n_J\le z^K<e^x$ and $J\ge z^{\del/2}/(2\log^2z)-1$,
we conclude that the lemma holds with $c=\del/3$, $a=499\del$ 
once $x$ is sufficiently large.
\end{proof}

%
%

\begin{lem}\label{populargood}
There is an absolute constant $c>0$, so that if
$\a>0$, $x$ is large (depending on $\a$) and $(\a,\frac1{10})$-good, then
there are $\gg \log x$ integers $n\le e^{x}$
satisfying $A(n) > n^c$ and $B(n)>n^c$.
\end{lem}

\begin{proof}
Let $\eps={1}/{10}$.  Let $\del$, $\gamma$, and $\eta$ be the constants from
Lemma \ref{E}.   Define $\TT$ as in \eqref{Tdef},
$\SS$ as in \eqref{S} and put $\widetilde{\SS}=\{p\in \SS : p \ge \sqrt{x} \}$.
Let $N:=\# \widetilde{\SS}$, so that from \eqref{Sb} we have 
$N\ge(\gamma/4)x/\log x$ for $x$ large.  Also, $N\le 2x/\log x$.
Let $\mathcal{Q}$ be the set of primes $q\le x^{1/2-\del}$ with $q\not\in\TT$.
For $q\in \mathcal{Q}$,
by \eqref{vqn} and the Brun--Titchmarsh inequality, we have
\be\label{Nq}
N_q := \# \{ p\in \widetilde{\SS}: q\mid p+1 \} \asymp\frac{N}{q}.
\ee
Suppose $k$ is an integer with $N^{1/2}\le k\le N^{3/4}$.
For $q\in\QQ$, if we choose a $k$-element subset $\MM$ of 
$\widetilde{\SS}$ at random,
we expect that the number of $p\in\MM$ with $q\mid p+1$ to
be $kN_q/N$.  That is, we are viewing a prime $p$ as corresponding to
the random variable which is 1 if $q\mid p+1$ and 0 otherwise.
By a standard result in the theory of large deviations
(see \cite[Sec.~5.11, (5)]{GS})
we have that the number of choices of $\MM$ with
\be\label{pMq}
\# \{ p\in \mathcal{M} : q\mid p+1 \} \ge \frac{k N_q}{2N}
\qquad \hbox{for all }q\in \mathcal{Q}
\ee
is at least, for some absolute positive constant $\nu$,
$$\left(1-\sum_{q\in\QQ}e^{-\nu kN_q/N}\right)\binom{N}{k}
\ge\frac12\binom{N}{k}\ge\frac12\left(\frac{N}{k}\right)^k
$$
for large $x$.  
(That the probabilistic model has us choosing
``with replacement" is easily seen to be negligible).
As in the proof of the previous lemma,
$n(\MM)=\sigma( \prod_{p\in \MM}p)<x^k$
and $P(n(\MM))\le x^{1/2-\del} < \log(x^k)$.
By \eqref{smooth}, there are $\le x^{k/30} \le N^{k/29}$
distinct values $n(\MM)$.  Hence, for large $x$ there is 
some integer $n<x^k$ with many representations as $n(\MM)$ where
$\MM$ satisfies \eqref{pMq}; in particular
$$
B(n) \ge \frac12 \pfrac{N}{k}^k N^{-k/29} \ge x^{k/5} > n^{1/5}.
$$

We next show that for each such $n$ we have $A(n)$ large.  
Note that generalizing \eqref{rad}, we have that if $w$ is
a positive integer with $\phi(w{\cdot}\rad(n))\mid n$, then
$$
n=\phi\left(w{\cdot}\rad(n)\frac{n}{\phi(w{\cdot}\rad(n))}\right).
$$
Thus, we can show that $A(n)$ is large if we can show that there
are many such integers $w$ with $(w,n)=1$
(to ensure that the integers 
$w{\cdot}\rad(n){\cdot}n/\phi(w{\cdot}\rad(n))$ are distinct for 
different $w$'s).  Towards this end, let
$$
\SS'=\{p\le x: p>\sqrt{x},~q\mid p-1\hbox{ implies }q\in\QQ\},\quad N'=\#\SS'.
$$
By Lemma \ref{E} and \eqref{sum1/t} we have $N'\gg x/\log x$,
so that $N'\asymp N$.
For each $q^j$ with $q\in\QQ$, let
$$
N_{q^j}':=
\#\{p\in{\SS'}:q^j\| p-1\}
$$
so that the Brun--Titchmarsh inequality implies that
$N_{q^j}'\ll x/(q^j\log(ex/q^j))$ for $q^j\le x$.
Put $k'=\lceil \xi k \rceil$, where $\xi$ is a small, fixed positive number.
For each $k'$-element subset $\MM'$ of ${\SS'}$, let 
$w(\MM')=\prod_{p\in\MM'}p$.  If $\MM'$ is chosen at random, 
the expected value of $\sum_{p\in\MM'}v_q(p-1)=v_q(\phi(w(\MM')))$ 
is $k'\sum_{j\ge1}jN_{q^j}'/N'$  (we are now viewing our random
variable as $v_q(p-1)$).
By the same result in \cite{GS},
there are at least $\frac12\binom{N'}{k'}$ choices for $\MM'$ with
$$
v_q(\phi(w(\MM')))\le\frac32k'\sum_{j \ge 1}\frac{jN_{q^j}'}{N'}
\qquad \hbox{ for all }q\in\QQ.
$$
For such choices of $\MM'$, we have $v_q(\phi(w(\MM')))\ll k'/q$, so if
we choose $\xi$ small enough, we have
$$
v_q(\phi(w(\MM')))\le k \frac{N_q}{4N}\le \frac12v_q(n),
$$
by \eqref{Nq} and \eqref{pMq}.  Since (cf. \eqref{vq})
$$
v_q(\phi(\rad(n)))\le \frac{x^{1/2-\del}}{q-1} \le \frac12 v_q(n),
$$
and since each prime factor of $w(\MM')$ is $>x^{1/2}\ge P(n)$, 
we deduce that
$\phi(w(\MM'){\cdot}\rad(n)) \mid n$ and that the numbers
$w(\MM'){\cdot}\rad(n){\cdot}n/\phi(w(\MM'){\cdot}\rad(n))$ 
are distinct for different choices of $\MM'$.  It follows that 
$$
A(n) \ge \frac12 \binom{N'}{k'} \ge \frac12 \pfrac{N'}{k'}^{k'}
> x^{k'/5} \ge n^{\xi / 5}.
$$
Put $c=\min(1/5,\xi/5)$.
Notice that our construction of $n$ depends on $k$, and
$$
x^{k/2} \le n \le x^k\le e^x.
$$
Letting $k$ run over the powers of 2 in $[N^{1/2},N^{3/4}]$ produces $\gg\log x$
distinct values of $n$, each $\le e^x$, for which $A(n)>n^c$ and $B(n)>n^c$.
\end{proof}

%
\section{Further problems}\label{sec:further}
%

\begin{enumerate}
\item
It is known that for any integer $k\ge 1$, there are integers $n$ with $B(n)=k$
and for any integer $l\ge 2$, there are integers $n$ with $A(n)=l$, see
\cite{F99}, \cite{FK}.  The famous Carmichael conjecture states that 
$A(n)$ is never 1, but this is still open. 
\begin{conj}
For every $k\ge 1$ and $l\ge 2$, there are integers $n$ with 
$A(n)=l$ and $B(n)=k$.
\end{conj}
\noindent
Schinzel has shown (private communication; see also \cite{S61}) 
that this conjecture follows from his Hypothesis H.

\item
If, as conjectured by Hardy and Littlewood, the number of pairs of
twin primes $\le x$ is $\sim C x/\log^2 x$, then the number of common values
$n\le x$ of $\phi$ and $\sigma$ is $\gg x/\log^2 x$.   
What is the correct order of
$\# \{ n\le x: A(n)\ge 1\text{ and }B(n)\ge 1\}$ ?

\item
Does $\phi(a)=\sigma(b)$ have infinitely many solutions with squarefree integers
$a,b$?  Our construction, when using $(\a,\eps)$-good values of $x$,
uses squarefree $b$ while $a$ is divisible by large
powers of primes.
\medskip

\item
As mentioned, Erd\H os showed that $A(k!)\ge 1$ for every positive integer 
$k$ \cite[sec. B38]{Guy}.  Is $B(k!) \ge 1$ for every $k\ne2$? 
How about at least infinitely often?  Note that our proof in 
Lemma~\ref{populargood} shows that there is some number $c>0$ such
that $A(k!)\ge (k!)^c$ for every $k$.
\end{enumerate}

\bigskip

{\bf Remarks.} 
There is an alternative approach to proving Theorems \ref{thm1} and \ref{thm2}
(with a somewhat weaker conclusion about the number of common values below $x$),
suggested to us by Sergei Konyagin.
Namely, it is possible to prove, using Lemmas \ref{Landau} and \ref{Page},
that there is an $\a>0$ such that for large $u$, there is a value of
$x\in [\log u,u]$ which is $(\a,\frac1{10})$-good.  Indeed,
let $\lam>0$ be small, and let $\a$ be the constant from Lemma \ref{Page}.
Let $\gamma$ be a constant satisfying $\gamma > 1/(10\a)$.
Let $m_1,m_2,\ldots$ be the (possibly empty) list of moduli for which
there is a character $\chi\in C(m_j)$ and zero
 $\beta_j\ge 1 - {\lam}/{\log m_j}$ of $L(s,\chi)$.  Let $j$ be the largest
index with $m_j \le (\log x)^{\a}$.  If there is no such $j$, then
$x$ is $(\a,\frac1{10})$-good.   Otherwise, $u=\max(\log x,\exp \{
\gamma (1-\beta_j)^{-1} \})$ is $(\a,\frac1{10})$-good upon using the
definition of $j$ and applying
Lemma \ref{Landau}.


\end{document}